\documentclass[a4paper,12pt]{amsart} \usepackage{euscript,eufrak,verbatim}
\usepackage[final]{epsfig}
\usepackage{psfrag}

\def\a{\alpha}

\def\Z{\mathbb{Z}} 
\def\R{\mathbb{R}}

\def\Ai{A^{(i)}}

\def\C{\mathcal{C}}
\def\V{\mathcal{V}}
\def\U{\mathcal{U}}
\def\M{\mathcal{M}}
\def\P{\mathcal{P}}

\newtheorem{theorem}{Theorem}
\newtheorem{lemma}[theorem]{Lemma}
\newtheorem{proposition}[theorem]{Proposition}
\newtheorem{corollary}[theorem]{Corollary}

\DeclareMathSymbol{\varnothing}{\mathord}{AMSb}{"3F} \begin{document}
\renewenvironment{proof}{\noindent {\bf Proof.}}{ \hfill\qed\\ }
\newenvironment{proofof}[1]{\noindent {\bf Proof of #1.}}{ \hfill\qed\\ }

\title{Periodic billiard orbits in right triangles}  
\author{Serge Troubetzkoy}
\address{Centre de physique th\'eorique\\
Federation de Recherches des Unites de Mathematique de Marseille\\
Institut de math\'ematiques de Luminy and\\ 
Universit\'e de la M\'editerran\'ee\\ 
Luminy, Case 907, F-13288 Marseille Cedex 9, France}
\email{troubetz@iml.univ-mrs.fr}
\urladdr{http://iml.univ-mrs.fr/{\lower.7ex\hbox{\~{}}}troubetz/} \date{}
\begin{abstract}  
There is an open set of right triangles such that for each
irrational  triangle
in this set (i) periodic 
billiards orbits are dense in the phase space, 
(ii) there is a unique nonsingular perpendicular billiard
orbit which is not periodic, and (iii) the perpendicular periodic orbits fill 
the corresponding invariant surface.
\end{abstract} 
\maketitle

\pagestyle{myheadings}

\markboth{Billiards in right triangles}{SERGE TROUBETZKOY}

\section{Introduction}

A billiard ball, i.e.~a point mass, moves inside a polygon $Q  \subset \R^2$ 
with unit speed along a straight line until it reaches the boundary
$\partial Q$, then instantaneously changes direction according to the mirror
law:  ``the angle of incidence is equal to the angle of reflection,'' and
continues along the new line. If the trajectory hits a corner of the polygon,
in general it does not have a unique continuation and thus by definition
it stops there.

Billiards in polygons are easy to describe, but it is difficult to prove
deep theorems about them because of a lack of machinery.
For example, it is unknown if every polygon
contains a periodic billiard orbit.  On the other hand, for so called
rational polygons, one can apply Teichm\"uller theory to obtain many
deep theorems.  It is known that all rational polygons possess many
periodic orbits \cite{Ma}, and in fact they are dense in the phase
space \cite{BoGaKrTr}.  Galperin, Stepin and Vorobets proved
many other interesting results about periodic orbits \cite{GSV}.
Many of these results can be
found in the introductory book by 
Tabachnikov \cite{T}, and in several
survey articles \cite{G1,G2,MT}.

The main result of this article produces
the first non rational polygons for which the above mentioned result is true: 
{\em periodic billiard trajectories are dense in the phase space
for an open set of right triangles} (Theorem \ref{thm0}). 
This result 
is proved by a careful analysis of the symmetries of perpendicular 
periodic orbits (Theorem \ref{thm1}).
Using these symmetries we prove that {\em there is a unique perpendicular escape
orbit} (Theorem \ref{thm2})
and that {\em the invariant surface which contains the perpendicular
direction is completely foliated by perpendicular periodic orbits}
(Theorem \ref{pascal}).
The billiard in a right triangle is  well known to be equivalent
to the mechanical system of two elastic point particles in an interval (see for
example \cite{T}), thus our results hold for this system as well.
 
\section{Statement of results}

Cipra, Hanson and Kolan have shown that almost every orbit which is 
perpendicular to the base of a right triangle is periodic \cite{CHK,T}.  
Here the
almost every is with respect to the length measure on the side of the
triangle in question.
Periodic billiard orbits always come in
{\em strips}, i.e.~for any $x = (q,v)$ whose billiard orbit is periodic
where $ q \in \partial Q$ and $v$ is any inward pointing direction
there is an open interval $I \subset \partial Q$ such that $q \in I$ and
for any $q' \in I$ the orbit of $x' = (q',v)$ visits the same sequence of sides
as $(q,v)$ and thus in particular is periodic (see Figure 1). 
\begin{figure}
\centerline{\psfig{file=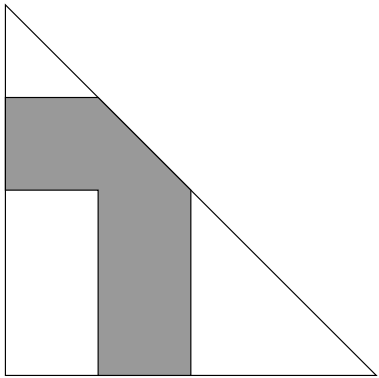,height=50mm}}
\caption{A periodic strip.}
\end{figure}  
A maximal width strip will be called a {\em beam}.
All orbits in a periodic beam
have the same period except perhaps the middle orbit which has half
the period in the case its period is odd.
The results of \cite{CHK} imply that a set of
full measure of the base of the triangle is covered by an at most countable
union of intervals such that each of these intervals forms a periodic beam.

They mention
``our computational evidence also suggest that the trajectory at the 
middle of each (perpendicular periodic) beam hits the right-angle vertex 
of the triangle.''  This had also been noticed in an earlier article
by Ruijgrok \cite{Ru}.
They also speculate that there is at most one nonsingular, non periodic 
trajectory. These speculations were the starting point of this research,
for this article. I present an elementary proof of these facts for 
irrational right triangles with smaller angle $\alpha$  satisfying 
$\frac{\pi}{6} < \alpha <\frac{\pi}{4}$.

Make the {\em convention} that one leg of the right triangle
is horizontal with the smaller angle $\alpha$ of
the right triangle being the angle between this leg and the hypotenuse, and  
that the word perpendicular (orbit, beam, etc.)
refers to perpendicularity to the horizontal leg. 

Any billiard trajectory which hits a right-angled vertex (or more generally
a vertex with angle $\pi/n$ for some positive integer $n$) has a unique
continuation.  Reflect the triangle in the sides of right angle to obtain
a rhombus.  The study of the billiard in the triangle reduces to that
in the rhombus (see next section for details).  
Throughout the article all beams considered will 
be with respect to the angle coding in the rhombus defined in the next 
section. 

\begin{theorem}\label{thm1}
For any irrational right triangle  whose smaller angle satisfies
$\frac{\pi}{6} < \alpha <\frac{\pi}{4}$
consider any perpendicular periodic beam of period $2p$. Then
\begin{enumerate}
\item{} the midpoint of the
beam hits the right-angle vertex of the
triangle (i.e.~the mid point of the rhombus), 
\item{} the beam returns to itself after half
its period with the opposite orientation, and
\item{} $p$ is even and the first $p+1$ letters of the code of the beam
form a palindrome.
\end{enumerate}
\end{theorem}

The theorem implies that when viewed as an object in the
phase space of the billiard flow in the triangle
the beam is a M\"oebius band.  This is not in contradiction with
the well known construction,
of invariant surfaces (see for example \cite{G1,T} in the
rational case and \cite{GT} in the irrational case) since
the directional billiard flow is isomorphic to the geodesic
flow on the invariant surface for any direction except 
for directions which are parallel to a side
of the polygon. The direction we are considering is such a direction.

The symmetry of the beam is reminiscent
of the symmetry of the beam of periodic orbits around a periodic
orbit of odd period mentioned above.
Galperin, Stepin and Vorobets showed that any perpendicular orbit whose
period is not a multiple of 4 is unstable under perturbation \cite{GSV}.  
Viewed as an orbit in the rhombus the period is a multiple of 4, however
a simple argument shows that the  period is not a multiple of 4 when
the orbit is viewed as an orbit in a right triangle \cite{GZ}.

Call an orbit {\em recurrent} if its code $\{a_i\}$ satisfies $a_j = a_0$
for some $j > 0$. Call a nonrecurrent orbit an {\em positive escape orbit} 
if  $\limsup _{i \to \infty} a_i = \infty$ and {\em negative escape orbit}
if  $\liminf _{i \to \infty} a_i = -\infty$.
Call a direction {\em simple} if there are no generalized
diagonals in the invariant surface containing the direction.
Here the invariant surface is thought of as arising from the right
triangle, so that the vertex at the right angle is considered as
a singularity.

\begin{theorem}\label{thm2}  $\phantom{1}$
\begin{enumerate}
\item{} Consider an irrational right triangle.  In any  simple direction 
then there is at most one non-singular positive escape orbit
and at most one non-singular negative escape orbit.
\item{}  For any irrational right triangle  whose smaller angle satisfies
$\frac{\pi}{6} < \alpha <\frac{\pi}{4}$,
there is at most one non-singular escape orbit in the perpendicular direction. 
If it exists it is both a positive and negative escape orbit.
\end{enumerate}
\end{theorem}

If we do not care if the orbit is singular or not, then we can remove the
words at most from the statement of the theorem.
Theorem \ref{thm2} is a substantial strengthening of a particular case 
of a result of \cite{GT}
who show that the set of nonrecurrent orbits has measure zero and
of a particular case of a result of \cite{ST} who show that this set has lower box 
counting dimension at most one half. 

Boshernitzan has conjectured that given a rational triangle, every nonsingular
orbit is periodic in  the invariant surface
containing a perpendicular direction \cite{Bo}.  I prove the following
analog for irrational right triangles.

\begin{theorem}\label{pascal}
Fix an irrational right triangle whose smaller angle satisfies
$\frac{\pi}{6} < \alpha <\frac{\pi}{4}$.
Consider the 
invariant surface of $M$ in the perpendicular direction.
Then all the nonsingular orbits on $M$
except the unique escape orbit are periodic.
\end{theorem}

Combining Theorem \ref{pascal}  with the main theorem of
 \cite{BoGaKrTr} yields

\begin{theorem}\label{thm0}
Periodic orbits are dense in the 
phase space of irrational right triangles whose smaller angle satisfies
$\frac{\pi}{6} < \alpha <\frac{\pi}{4}$.
\end{theorem}

\section{Definitions and proofs}

There is a nice introductory book on billiards by Tabachnikov \cite{T} and
several survey articles \cite{G1,G2,MT} which can
be consulted for details on polygonal billiards in general.
For a rational polygon there is a well known construction of
invariant surfaces, they are always compact.  They same construction
leads to noncompact invariant surfaces for irrational polygons.
The invariant surfaces for an irrational right triangle can be
thought of as having a $\Z$-quasiperiodicity \cite{CHK,GT,Tr,ST,T}.
 
\begin{figure}
\centerline{\psfig{file=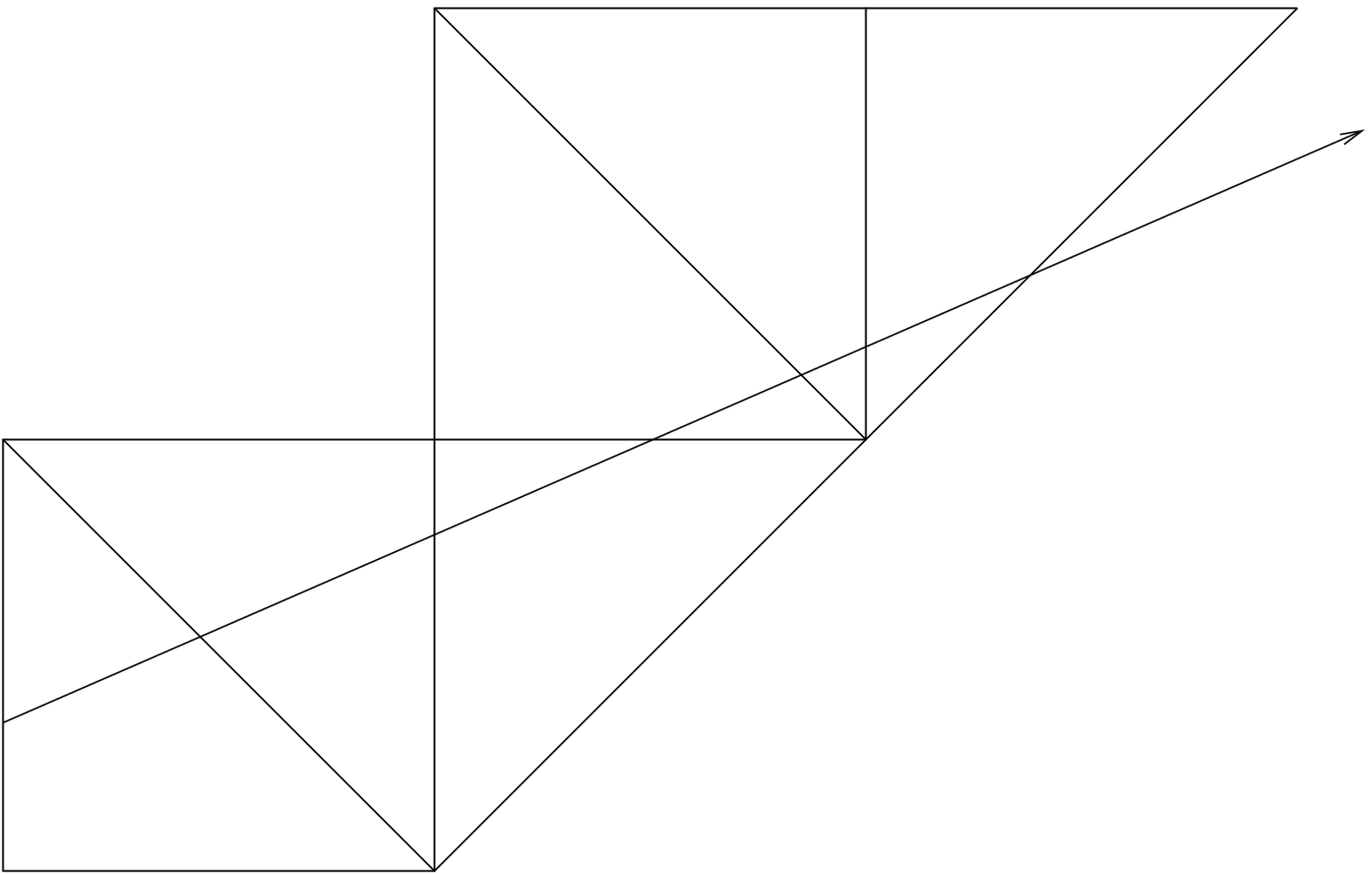,height=30mm} \hspace{.51in}
\psfig{file=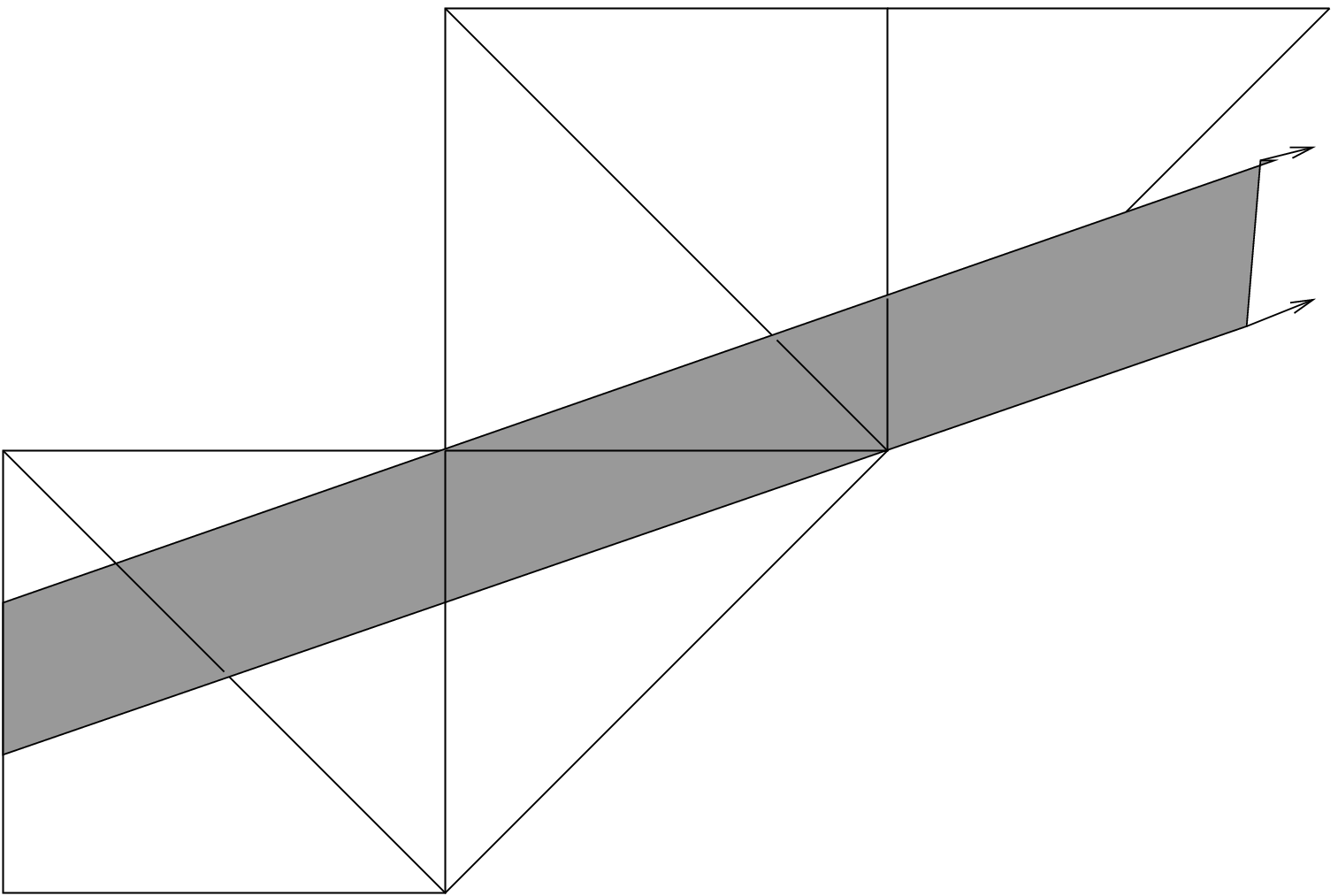,height=30mm}}
\caption{Unfolding a billiard trajectory and the associated strip.}
\end{figure}  

The proof is based on the procedure of {\em unfolding of a 
billiard trajectory}. 
Instead of reflecting the trajectory with respect to a side of a
polygon reflect the polygon in this side.  Thus the trajectory
is straightened to a line with a number of isometric copies of the 
polygon skewered on it (Figure 2a). 

Fix an orbit segment.
There is a {\em strip} around this trajectory segment 
such that the same sequence 
of reflections is made by all trajectories in the strip (Figure 2b). 
The number of reflections is called the {\em length} of the strip. 
Call a maximal width strip a {\em beam}.
The boundary of a periodic beam consists of one or more
trajectory segments which hit a vertex of the polygon. 
If this vertex is the right-angle one, then
the sequence of reflections on both sides of it is essentially the same 
since the singularity due to such a vertex is removable.  
This enables us to
consider the billiard inside a rhombus which consists of four copies
of the right triangle and unfold the rhombus (Figure 3).
The central and
reflectional symmetries of the rhombus play an essential part
in the proofs of the theorems.

\begin{figure}
\centerline{\psfig{file=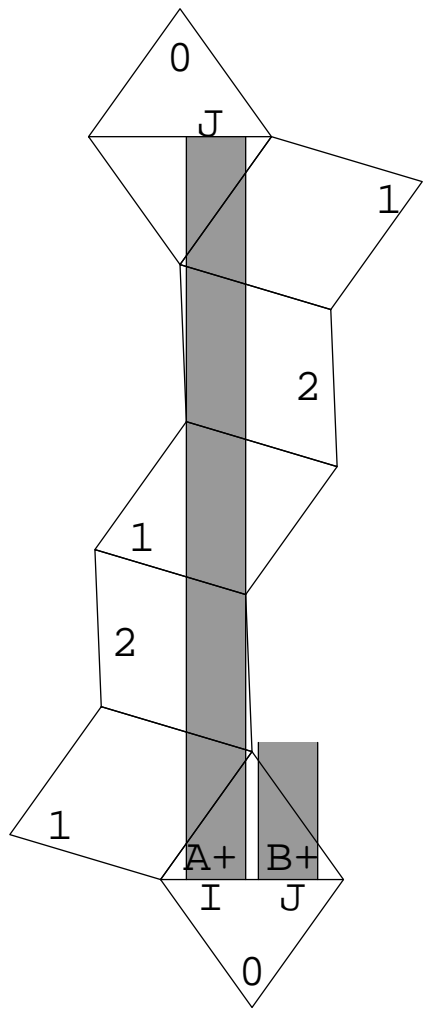,height=140mm}}
\caption{A periodic beam in a rhombus with code word $0121210$.}
\end{figure}  

\begin{figure}
\psfrag*{I}{\footnotesize$I_k$}
\centerline{\psfig{file=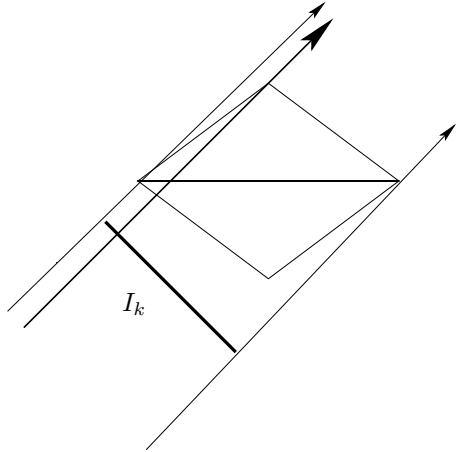,height=60mm}}
\caption{Rhombus $k$ projectivized to the interval $I_k$.}
\label{split}
\end{figure}

Throughout the article unfoldings will be of the billiard in the rhombus
and thus billiard orbits through its center will be considered as defined.
As we follow the reflections along a straight line trajectory
we see that each flip rotates the rhombus by $2\alpha$ where $\alpha$ is
one of the interior angles of the rhombus.  Thus we can {\em label} the rhombi
with integers according to the total number of (clockwise) rotations by
$2\a$. Each orbit is {\em coded} by the sequence of labelled rhombi it 
visits, for example any orbit in the beam $A^+$ in Figure 3 has code 0121210. 
It is often convenient to consider labelled rhombus in a projective
sense as an interval, namely for rhombus $k$ consider the beam of length 1
with code $k$ and then the associated interval $I_k$ is just the perpendicular
width of of the beam (see Figure \ref{split}). 
The $k$th-rhombus or the associated interval $I_k$ will often be
referred to as {\em level} $k$. 
All reference to the length of a beam will pertain
to the number of rhombi it crosses
which differs from its length in the right triangle.

Let us begin with the following fact.

\begin{lemma}\label{lemma1}
Suppose the $A^+$ and $B^+$ are parallel beams 
(not necessarily perpendicular or periodic)
with the same codes, then $A^+ = B^+$.
\end{lemma}

\begin{proof}
The proof is by induction.  Fix an initial direction.
The base case is clear, there is a single
beam with code $0$. Furthermore, there is a single beam with code $01$, 
a single beam with code $0(^{\scriptstyle -}\! 1)$ and these two beams of length two are
separated by an orbit which hits a vertex of the rhombus. 
Assume that every beam with length $n$ is the unique beam with its
code. Fix a beam of length $n$ and denote by $j$ the last rhombus it visits.
The $j$th rhombus is attached to the $(j+1)$st by two
parallel sides and to the $(j-1)$st
along the other two parallel sides. There are two points which
are common to the $(j-1)$st, $j$th and $(j+1)$st rhombi. During a given
pass through the $j$th rhombus, the
beam can hit only one of the two in exiting the $j$th rhombus (it
can hit the other in entering the $j$ rhombus).
The orbit which splits a beam is emphasized in Figure \ref{split}. 
If the beam does not arrive at this point
then it does not split and the beam can be continued to a unique beam of length
$n+1$. 
If the beam reaches this point then its continuation is split into 
two sub-beams of length $n+1$ with the resulting sub-beams code
differing in the $(n+1)$st place. 
\end{proof}

Next we need the following fact.

\begin{lemma}\label{lemma2}
The three statements in Theorem \ref{thm1} are equivalent.  
\end{lemma}

\begin{proof} Perpendicular periodic beams must 
have even period since they must retrace their orbits between
the two perpendicular hits.
Fix a perpendicular periodic beam of period $2p.$ 
Consider the beam  of length $p$
between the two perpendicular collisions (see Figure 3).

$(1 \Rightarrow 2$ and $3)$ 
Suppose that the midpoint of this beam hits the right-angle
vertex. The central symmetry of the rhombus at this collision point
implies that the beam itself is centrally symmetric in this point and thus
has palindromic code. Furthermore $p$ is even since the symmetry is 
about an interior point of the rhombus.
Label by $I$ the interval of the departure and by $J$ the interval of 
the arrival.
Both intervals are contained in the rhombus labelled 0, the
central symmetry of the beam implies that 
$I$ and $J$ are centrally symmetric.

$(2 \Rightarrow 1)$ By assumption the initial interval $I$ of the beam and
the final interval $J$  fold to the same interval in the triangle
Thus the midpoint of $I$ must have period $p$ or have a singular orbit 
as an orbit in the triangle.  
It is easy to see that the orbit can not be periodic, since any perpendicular
periodic orbit must be twice perpendicular and the collision at time $p/2$ is not
perpendicular. Thus the orbit is singular in
the triangle. Since it was nonsingular in the rhombus it hit the right-angled
vertex of the triangle.


$(3 \Rightarrow 1)$ Consider the
central symmetry of the palindromic beam around its center $c$. This takes the beam
to a priori another beam with the same code.  Lemma \ref{lemma1} implies
that the beams must coincide. Thus the beam is centrally symmetric about
its center. Since $p$ is even the point $c$
must lie in the interior of the rhombus visited by the strip at step $p/2$.
However, viewed in the rhombus
one sees that the beam can only be centrally symmetric around the center
of the rhombus and no other point, thus $c$ must be the
center of the rhombus.
\end{proof}

%
%
%
%
We now turn to a detailed analysis of beams which start in a given
direction, and either return to that direction or another fixed
direction.
\begin{lemma}\label{lemma5'}
Fix an angle $\theta_0$ and consider the angle coding of the billiard
with $\theta_0$ taken as level $0$.
Suppose $N > 0$.
Consider the set $S_0$ of points whose
beam of length 2 has code $01$.
This set is an interval, partition it  
into intervals such that for each partition interval there exists 
$p \ge 1$ such that the initial orbit code $\{a_i\}_{i=0}^p$ 
of all points in the interval does not depend on the points and 
is such that 
$a_0 = 0$, $0 < a_i < N$ 
for $i = 1, \dots p-1$ and $a_p$ is 
equal to either $0$ or $N$.
Then there are at most $|N|$ such codes and corresponding subintervals of $S_0$.
The interiors of these subintervals are
pairwise disjoint and their union covers $S_0$.
\end{lemma}

\begin{proof}
This lemma has been essentially proven in \cite{GT} and somewhat more
explicitly in \cite{ST}.
I repeat the proof here for completeness.   
Consider the beam of length two whose initial code is 
$01$.  Follow all the orbits
in this beam until they reach level $N$ or return to level $0$ at some time 
$p$.  For the beam to split, an orbit in the beam must reach a 
vertex.

This can be seen for example in Figure 3, the beam (not drawn in the picture)
just to the left of the beam $A^+$ has code $01210$ 
(the copy of the rhombus corresponding
to last $0$ is not drawn in the picture).  These two beams are split by a vertex
in (the $2$nd occurence) of the $1$st rhombus.

For each $i$ there is a single vertex $v_i$ which can split the beam in 
rhombus $i$. Since the billiard is invertible this vertex can only be 
reached from at most one point in the $0$th rhombus (without returning 
to the rhombus). The set
$\V := \V(N)  := \{v_i: 1 \le i \le N-1\}$ of vertices has cardinality $N-1$, 
therefore there are at most $N$ distinct beams and thus codes.
\end{proof}

\begin{lemma}\label{lemma3}
Fix an angle $\theta_0$ and consider the angle coding of the billiard
with $\theta_0$ taken as level $0$.
For any $M<0<N$ we can partition level $0$ into a finite number of intervals
with each interval in one of three classes, 
$\P:=\P(\theta_0,M,N),\M:=\M(\theta_0,M,N),\U:=\U(\theta_0,M,n)$ such that
the orbits of any point in an interval in 
class $\P$ and $\M$ does not reach level $M$ or $N$. Furthermore
\begin{enumerate}
\item{} the forward orbit of each boundary point of each interval is singular,
\item{} the interior of each of the intervals in $\P$ consists periodic points
and the boundaries of each periodic family consists of generalized diagonals,
\item{} the intervals in $\M$ consists of one (or several) minimal
interval exchange transformations (i.e.~on a finite number of intervals!)
and the boundary of each minimal component consist of generalized diagonals and
\item{} the orbit of each point in the interior of an interval in $\U$ reaches level $M$ or
$N$.
\end{enumerate}
\end{lemma}

We have the following immediate corollary.

\begin{corollary}\label{cor1}
The orbit of every point in level $0$ which does not reach level $M$ or $N$ is either
singular or returns to level $0$.
\end{corollary}

The classes may well be empty, for example if we choose a direction
$\theta_0$ for
which the surface has no generalized diagonals then both $\P$ and $\M$ must
be empty. Galperin has shown that for almost all right triangles
there exists a $\theta_0$ such that the set $\M$ is nonempty \cite{Ga}.

\begin{proof}
Consider the cover of $S_0$ defined in the previous lemma and use the
subintervals which return to define a partially defined IET on level $0$.
Do the same for points whose code starts with $0(-1)$.
Call the domain of definition $D$. 
Do the same for the backwards dynamics and remark that these two
partially defined maps are inverses of each other where ever they are defined. 
The total length of the intervals where the forward map is not defined is equal to the
total length of the intervals where the backwards map is not define.
Thus we can
(in an arbitrary manner) complete the definition of the partially
defined map to an IET,  which we call the ghost map. 
The ghost IET agrees with the
partially defined first return billiard map whenever it was defined.  

Since the ghost map is an IET with interval of definition being level $0$.
The well known topological decomposition holds,
the interval of definition is decomposed into periodic and minimal components,
with the boundary of the components consisting of generalized diagonals
(see for example \cite{KaHa}).

Any point whose orbit enters an interval on which the ghost dynamics differs
from the true dynamics will
reach level $M$ or $N$  following the
true dynamics. On the other hand the topological decomposition mentioned above
can be applied to the  point for which the ghost and true dynamics always agree.
We remark that a ghost minimal component will, by minimality be either completely
contained in $D$ or every orbit will enter an interval on which the ghost dynamics 
differs from the true dynamics and thus will
reach level $M$ or $N$  following the true dynamics. 
\end{proof}

We next begin the analysis of simple directions.  The main reason to analyze such
directions is as a ``warm-up'' for perpendicular directions.

\begin{lemma}\label{lemma4}
Suppose $\theta_0$ is simple. Then for all integers $M < N$ 
\begin{enumerate}
\item{} there exists a nonsingular orbit segment starting in level $M$ and ending in level $N$ and 
\item{} a nonsingular orbit segment starting in level $N$ and ending in level $M$.
\end{enumerate}
\end{lemma}

\begin{proof}
We first remark that if there is a nonsingular orbit from $M$ to $N$ then since the map is
locally an isometry there is a whole interval which maps from $M$ to $N$. Since a.e.~point
in level $M$ is recurrent this implies that there must be orbits from $N$ to $M$ as well.
Thus it suffices to show that either (1) or (2) holds.

Suppose that neither (1) nor (2) holds.  Consider the set $S_M$ of all points 
whose code starts with
$M(M+1)$ and $S_N$ whose code starts with $N(N-1)$.  
Let $B_M^N$ be the union of the levels $M+1$ to $N-1$.
Since $\theta_0$ is simple the previous lemma implies that the backwards 
orbit of every point in $B_M^N$  
must reach $S_N$, $S_M$ or be singular.  Thus we can partition $B_M^N$
into intervals whose backwards orbit hit $S_M$, intervals whose backwards orbit
hits level $S_N$, and the backwards singular orbits which divide these intervals. 

Since we have assumed there are no orbits connecting levels $M$ and $N$, Corollary
\ref{cor1} implies that if
the backwards orbit hits $S_M$ (resp.~$S_N$) then the forward orbit hits level $M$ (resp.~level $N$). 
This implies that the points which divide the partition of $B_M^N$  
are also forward singular. Since
they are both backwards and forward singular they are part of a generalized diagonal, which
contradicts that fact that $\theta_0$ is simple.  Therefore at least one of (1) or (2) holds, and thus both hold.
\end{proof}

Next we prove the following strengthening of Lemma~\ref{lemma5'} for simple directions.
This is the first place that the role of the center of the rhombus becomes apparent.

\begin{proposition}\label{lemma5}
Suppose $\theta_0$ is  simple. Then 
\begin{enumerate}
\item{} there is exactly one code for which $a_0 =0$ and $a_p = N$ and
\item{} for all codes for which $a_0 = a_p = 0$ the code is a palindrome and 
$p$ is even. 
\end{enumerate}
\end{proposition}

\begin{proof}
Consider the set $S_0$ of all points whose code starts with $01$ and  
the set $S_N$ of all points whose code starts 
with $N(N-1)$. Set $S := S_0 \cup S_N$. 
Consider the points in $S$ whose orbits arrive at a vertex in $\V$ before
returning to level $0$ or $N$. The cardinality of $\V$ is $N-1$, thus since
$\theta_0$ is simple we can apply Lemma \ref{lemma3} to conclude that
there are exactly $N-1$ such points.  They partition $S$ into $N+1$
intervals.  Call the associated
beams $\Ai$ where $1\le  i \le N+1$. 
Suppose the $i$th beam is of length $p(i)$, we denote 
its code by $\{a^{(i)}_j\}_{j=0}^{p(i)}$. 

Consider the set $\C := \C(N) := \{c_i: 1 \le i \le N-1\}$ where  $c_i$ is the center 
of the $i$th rhombus. Consider the orbits
which start in $S$ and arrive at one of these centers before 
returning to level $0$ or $N$. 
By simplicity there are $N-1$ such orbits and they
are disjoint from the orbits which start in $S$ and arrive at a vertex in $\V$.
The orbit of $c_i \in \C$ is symmetric, thus
the code of a beam $\Ai$ containing such a point
is a palindrome and $p$ is even.
The code begins and ends with the same symbol, $0$ or $N$, 
and  is strictly positive in between, therefore it can not be a palindrome twice. 
Thus each of the orbits which arrive at a center are in distinct beams.
This implies that  each interval contains at most one $c_i$ and
all but two of the beams contain a center.

We have shown that there are two intervals 
whose beam does not contain a center. All the other beams start and end on one of the two
levels $0,N$ without visiting the other. Thus by Lemma \ref{lemma4}, the two beams must 
be as follows:
one connects levels $0$ to $N$ and the other levels $N$ to $0$ and all other codes are palindromes.
\end{proof}

We would like to extend Proposition~\ref{lemma5} to non simple directions, especially to perpendicular
directions.  There are three places where simplicity was used in the proof, the first
is that the pull back of all $v_i \in \V$ and all $c_i \in \C$ reaches 
the set $S$.  A finer counting argument overcomes this difference. 
The second is that the pull
backs of the $v_i$ and $c_i$ are disjoint.
Finally we applied Lemma \ref{lemma4}.  We are able to avoid the last two difficulties
by assuming that the angle $\alpha$ satisfies $\alpha \in (\frac{\pi}{6},\frac{\pi}{4})$.

\begin{lemma}\label{babypascal}
Fix an irrational right triangle $Q$ whose smaller angle 
$\alpha \in (\frac{\pi}{6},\frac{\pi}{4})$.
Consider the angular coding with respect to the 
perpendicular direction.  Then the first return map to level 0 can 
not be an IET.
\end{lemma}

\begin{proof}
Let $M$ be the invariant surface containing the perpendicular direction.
Consider the set $G$ of generalized diagonals which are contained in
the invariant surface $M$. Note that since
$Q$ is irrational the set $G$ can be infinite.  Each generalized
diagonal is a segment in $M$.  Several generalized diagonals can 
intersect at a singular point, thus the set $G$ has the structure of a graph,
not necessarily connected.   

Suppose the conclusion is not true. Then, the
set of perpendicular periodic orbits form a finite union of
annuli $A \subset M$.
We have $A \ne M$ since the area of $M$ is infinite and the area of $A$ is finite.
Thus $A$ must have a boundary which must consist of a finite union of generalized 
diagonals.  
We will analyze $G$ to show that this boundary must be empty. 

The boundary of $A$ consists of periodic loops.  
A {\em periodic loop} is a finite union of generalized diagonals which is 
periodic. We can think of a periodic loop as a (one-sided) boundary  of
a family of periodic orbits. The idea of the proof is to analyze all possible 
periodic loops which could possibly form the boundary of $A$, i.e.~those that
have perpendicular periodic orbits on at least one side.  Our analysis
will show that in each case they must have perpendicular periodic orbits on 
both sides and thus they can not form part of the boundary of $A$.

A periodic loop $g$ is called {\em simple} if it consists of a single 
generalized
diagonal, or equivalently the connected component of $G$ containing
$g$ is $\{g\}$.  Let $L$ be the horizontal diagonal of the rhombus. 

\begin{figure}
\psfrag*{L}{{$L$}}
\psfrag*{a}{\footnotesize$\alpha$}
\psfrag*{0}{\footnotesize{$0$}}
\psfrag*{1}{\footnotesize{$1$}}
\psfrag*{g}{{$g_L$}}
\centerline{\psfig{file=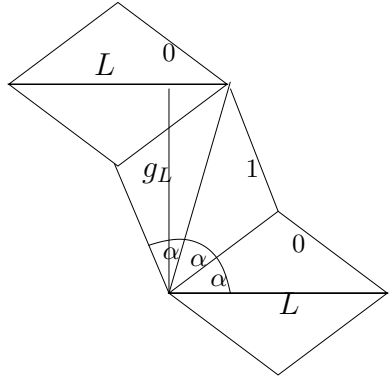,height=50mm}}
\caption{$g_L$ is a simple periodic loop.}
\label{fig5}
\end{figure}  
Consider the orbit $g_L$ starting from the left endpoint of $L$.  
Since $\pi/2 < 3\alpha < 3\pi/4$  the symbolic coding of $g_L$ is $010$,
and $g_L$  $g_L$ hits $L$ again in an interior point
without hitting any other vertex before (see Figure \ref{fig5}).  
Thus $g_L$ is a simple periodic loop.
Since $g_L$ hits $L$ at an interior point it has perpendicular periodic orbits on both 
sides of it and can not form part of the boundary of $A$.  
Similarly the perpendicular diagonal
starting at the right end point of $L$ can not be part of the boundary of $A$.

Consider any connected component $G'$ of $G$ such that 
$G' \cap \bar{A} \ne \emptyset$. 
This implies that there is at least one generalized diagonal $g_1 \in G'$ whose orbit
hits $L$ perpendicularly, and thus  we have a perpendicular periodic loop
$g'=g_1,\dots,g_n$ in $G'$.
Since we have finished treating it, we can
assume $g_L \not \in G'$ and thus $g_1$ hits $L$ in an interior point.

I claim that the periodic loop hits $L$ again perpendicularly at another interior
point of $L$. 
First of all $g'$ can not hit $L$ again at an endpoint by assumption.
Start flowing $g' (g_1)$ forwards and backwards
from $L$, this is a symmetric orbit in $M$.
At exactly half of $g'$'s (minimal) period, the forward and backwards orbits 
meet.  The meeting point must also be a symmetry point of the orbit.  The only points
of symmetry available are perpendicular reflections (in our case since $Q$ is irrational
the only axis of reflection in $M$ is $L$) and central symmetries through copies of
the center of $Q$. I claim that the only possibility for the second symmetry of $g$ is 
that it must hit $L$ a second time at an different interior point of $L$. 
First it can not hit a central
symmetry point since such points are isolated, while the second symmetry holds for the
whole periodic family which $g'$ bounds.  Furthermore $g'$ 
can not hit $L$ a second time at the same end point, in fact assuming this holds 
halves the period of $g'$, a contradiction of the definition of minimal period.  

\begin{figure}
\psfrag*{L}{$L$}
\psfrag*{g_1}{$\ g_1$}
\psfrag*{g_2}{$g_2$}
\psfrag*{g_3}{$g_3\ $}
\centerline{\psfig{file=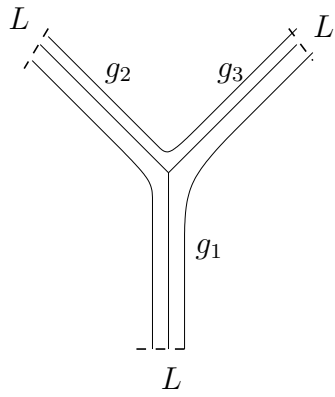,height=50mm}}
\caption{$G'$ consisting of three periodic loops, $\{g_1,g_2\}$, $\{g_2,g_3\}$,  and 
$\{g_3,g_1\}$. }
\label{fig6}
\end{figure} 

Thus both sides of $g_1$ must consist of periodic perpendicular orbits (see 
Figure \ref{fig6}).  If $g'$ is
simple then both sides of $g'$ consist of periodic perpendicular orbits and $g'$
can not be part of the boundary of $A$. If $g'$ 
is not simple then the ``other'' side of $g'$ is also a periodic loop, for which,
by the same reasoning, bounds perpendicular period orbits.  Since $A$ is a finite
union of annuli, there can only be a finite number of perpendicular periodic loops.
Thus, by
exhaustion, all of $G'$ consists of perpendicular periodic loops, and can not
form part of the boundary of $A$.  
This contradiction completes the proof.
\end{proof}

Next we prove the analog of Proposition \ref{lemma5} for perpendicular
directions.
\begin{proposition}\label{lemma12}
Fix an irrational right triangle $Q$ whose smaller angle 
$\alpha \in (\frac{\pi}{6},\frac{\pi}{4})$.
Consider the angular coding with respect to the 
perpendicular direction, then
\begin{enumerate}
\item{}  there is exactly one code for which  $a_0 = 0$ and $a_p = N$ and
\item{} for all codes for which $a_0 = a_p = 0$ the code is a palindrome and 
$p$ is even. 
\end{enumerate}\end{proposition}

\begin{proof}
We consider the same notation as in the proof of Proposition~\ref{lemma5},
the difference being that generalized diagonals exist. 
Suppose for the moment that no generalized 
diagonal passes through both $\C$ and $\V$.
Note that by symmetry a perpendicular orbit can reach only a single $c \in \C$.
Let $J \le N-1$ be the
cardinality of the points in $S$ whose orbits arrive at a vertex in $\V$ before
returning to level $0$ or $N$. They partition the $S$
into $J+2$ intervals and corresponding beams $\Ai$.  Let $K$ be the 
cardinality of the set of orbits which  start in $S$ and arrive $\C$.
Consider the beam $\Ai$ containing such an orbit, arguing as in the simple case 
that each center is in a different beam yields $K \le J+2$.

If no orbit from $S$ reaches $c_i$ before returning level $0$ or $N$,
then the orbit of $c_i$ never reaches rhombus $0$ or $N$.
Thus by Lemma \ref{lemma3}
it's orbit is periodic or its orbit closure is a minimal IET.
In either case since the map is a local isometry
we can associate $c_i$ with a vertex $v_{j(i)} \in \V$
such that no orbit from $S$ reaches this vertex.
Furthermore one generalized diagonal can not be the complete boundary to 
two minimal or periodic components of an IET.  Thus for a set of 
cardinality $L$ of such $c_i$ we must have at least $L$ distinct
such vertices $v_{j(i)}$.  This implies $N-1-J \ge  N-1-K$, or $J \le K$.
We have shown that $J \le K \le J+2$, i.e.~the number of beams which do not 
contain a center is at most 2.
These two beams are the only candidates for beams connecting level $0$ to $N$
and $N$ to $0$.

Now we deal with the case when there is an orbit passing through 
both $\C$ and $\V$.
Note that any orbit which starts perpendicularly and arrives at $C$ before $\V$,
returns, by the central symmetry to a perpendicular collision and is thus simple.
Thus we have reduced to the case when an orbit visits $\V$ before $\C$. 
By the central
symmetry it is a generalized diagonal. The counting
argument above works with the following differences. First of all 
we consider this generalized diagonal as a degenerate beam $\Ai$.
For this degenerate beam we interpret $J$ and $K$ as follows.
Since $\Ai$ visits a single $c \in \C$ it contributes one to $K$ and since
$J+2$ is the number of beams it contributes one to $J$.
With these modifications the inequality $J \le K \le J+2$ remains true
and  as above we conclude that
all but at most two of the beams contain a center.
By the central symmetry around the center the code of each such beam must 
satisfy
$a_0 = a_p \in \{0,N\}$. 

Next we need to show that exactly two exceptional beams exist. 
If no exceptional beam exists
then Lemma~\ref{babypascal} yields an immediate contradiction. If one 
exceptional beam exists, if $a_0=a_p$ then again
Lemma~\ref{babypascal} yields an immediate contradiction.
In the other cases ($a_0=0, \ a_p=N$ or $a_0=N, \ a_p=0$) 
the points of the exceptional beam
can not be recurrent, which is
a contradiction of the fact that almost every point is recurrent.
Thus exactly two exceptional beams exist.
All the other beams contain a center and thus are palindromes with
$p$ even.
 
Finally 
we need to show that in fact one of the two beams satisfies 
$a^1_0 = 0, a^1_p = N$
and the other satisfies $a^2_0 = N, a^2_p = 0$.
If we think of the beams as marked (1st and 2nd) then there are 
16 possible cases depending on 
of $a_i^j = 0$ or $N$ for $j=1,2$ and $i=0,p$.  Two of these cases are
good, we must eliminate all the 14 other cases.

Since almost every point is recurrent we can 
not have a beam
starting with $0$ and ending with $N$ without the other beam connecting level $N$
to level $0$, this eliminates 6 cases.

To be able to apply Lemma~\ref{babypascal} repeat the construction for points which
are perpendicular to the leg of the triangle but whose code is strictly negative
before returning to level $0$ or $-N$ to produce two addition exceptional beams 
with codes $b^j_i$,
which are centrally symmetric to the two original ones in the sense that
$a^j_i = -b^j_i$ for $j=1,2$.
A direct application of Lemma~\ref{babypascal}
eliminates the the remaining 8 bad cases when there is no beam starting in $0$ 
ending in $N$.
%
\end{proof}

%
%
%
%
%
%
%
%

\begin{proofof}{Theorem \ref{thm1}}
The theorem follows immediately for by combining 
Lemma~\ref{lemma2} and Propostion~\ref{lemma12}.
\end{proofof}

\begin{proofof}{Theorem \ref{thm2}} We start with the simple case.
First consider orbits whose code starts with $01$. Any such
nonsingular nonrecurrent orbit must be contained in the intersection 
$\cap_{N \ge 1} B^{i_0(N)}$ where $B^{i_0(N)}$ is 
the unique interval which satisfy alternative (1) for fixed $N \ge 1$  of 
Proposition~\ref{lemma5}.
Consider the initial segment of these beams, these intervals are clearly 
nested,  thus there is at most one interval or a single point in the intersection 
of their closures.   Since the set of nonrecurrent orbits has measure 
zero it consists of at most a single point. The orbits of
this point is either a positive escape orbits or singular.
Similarly we find a single positive escape orbit starting on any level $N$.
This escape orbit must, by uniqueness, coincide with the one starting on level $0$.

A symmetric argument holds for the backwards escape orbit, by starting with points
whose code commences with $0-1$.   

We turn to the perpendicular case.  The fist difference to the simple case is that
we replace Proposition~\ref{lemma5} by Proposition \ref{lemma12}. We immediately obtain that
there is at most one positive escape orbit starting on level $0$.  To show that this is
the unique positive escape orbit on the whole surface we argue like in the proof 
of Propostion~\ref{lemma12}.  Fix $M > N > 0$.  
The first part proof goes through for levels $M$ and $N$,
i.e.~we know that there are
at most two exceptional beams.  We need to conclude that they exist and that
they connect levels $M$ and $N$
to each other without applying Lemma~\ref{babypascal}. This follows immediately from
the fact that there is a beam connecting level $0$ to $N$ and a beam connecting level
$N$ to level $0$.

By symmetry, the positive and negative escape orbits must
be the forward and backwards part of the same orbit.  
\end{proofof}

\begin{proofof}{Theorem \ref{pascal}}
The  proof is similar to the proof of Lemma~\ref{babypascal}.
The set of perpendicular periodic orbits form a countable union of
annuli $A \subset M$.
If $A \ne M$ then $A$ must have a boundary which must consist of generalized 
diagonals and escape orbits. Besides the escape orbit, there is an additional
difficulty that the relevant part of the graph of generalized diagonals on 
$M$ need not be finite.

The argument that a {\em finite} connected component $G'$
satisfying  $G' \cap \bar{A}  \ne \emptyset$ 
can not form part of the boundary of
$A$ is the same as in the proof of Lemma~\ref{babypascal}.
This completes the proof in the rational case.

Next suppose that the escape orbit exists, i.e.~is nonsingular.
Denote this orbit by $e$.  By construction,
the orbit $e$ starts on level $0$ at a point $e_0$ which is an interior point of $L$. 
Thus we can find a sequence of perpendicular periodic orbits 
which approaches $e$ from the left and a sequence which approaches
$e$ from the right, and $e$ can not form part of the boundary of $A$.

\begin{figure}
\psfrag*{L}{$L$}
\psfrag*{g1}{$g_1$}
\psfrag*{g2}{$g_2$}
\psfrag*{g3}{$g_3$}
\psfrag*{e}{$e$}
\centerline{\psfig{file=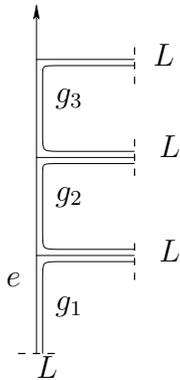,height=50mm}}
\caption{An infinite connected component $G'$ which includes the escape orbit. }\label{fig7}
\end{figure} 

The remaining case is $G'$  infinite.  There are two subcases, either every
edge in the graph has periodic perpendicular orbits on both sides or not. 
In the first case we argue as in the proof of Lemma~\ref{babypascal}, using
induction instead of exhaustion to complete the proof.  
In the second case the (unique) escape orbit $e$
is infact singular and contained in $G'$, see Figure \ref{fig7}.  
Thus the one side of $e$ bounds perpendicular
periodic families, while the other side is a limit  or perpendicular periodic families.
\end{proofof}

The method of proof also allows to verify the Boshernitzan conjecture for an infinite
collection of rational right triangles, 
those where the perpendicular orbits starting at the (non-right angle) verteces
are simple periodic loops. 

\begin{proofof}{Theorem \ref{thm0}}
For irrational triangles this follows immediately from Theorem~\ref{pascal} since
the surface $M$ is dense in the full phase space.  For rational triangles the
result was established in \cite{BoGaKrTr}.

\end{proofof}

\section{Extensions}
I would like to remark on the assumption $\frac{\pi}{6} < \alpha <\frac{\pi}{4}$.
This is a technical assumption which guarantees that the orbits starting at the
endpoints of $L$ are simple periodic loops.  This is the key
condition which is used in the proof of the theorems. 
It seems likely this condition holds for an open set of
full measure of parameters $\alpha$, but to prove this seems 
messy and complicated.

\section{Acknowledgements}
I profitted greatly from discussions with Pascal Hubert.

\end{document}